\RequirePackage{fix-cm}
\documentclass[smallextended]{svjour3}       
\smartqed  
\usepackage{amsfonts}
\usepackage{xr}
\usepackage{color}
\usepackage{graphicx,epsfig}
\usepackage{amsmath}
\usepackage{color}
\usepackage{hyperref}
\usepackage{tcolorbox}

\usepackage{tikz}
\usetikzlibrary{patterns}
\usepackage{float,subcaption}
\usepackage{lipsum}
\usetikzlibrary{decorations.pathmorphing} 
\usetikzlibrary{matrix} 
\usetikzlibrary{arrows} 
\usetikzlibrary{calc} 
\usepackage{cancel}
\usepackage[framemethod=TikZ]{mdframed}

\usepackage{natbib}

\newcommand{\R}{\mathbb{R}}

\newtheorem{lem}{\bf Lemma}[section]
\newtheorem{prop}[lem]{\bf Proposition}

\begin{document}

\title{A Turing mechanism in order to explain the patchy nature of Crohn's disease}


\author{ Gr\'egoire Nadin  \and Eric Ogier-Denis \and Ana I. Toledo  \and  Hatem Zaag }


\institute{G. Nadin \at
              Laboratoire Jaques-Louis Lions, Universit\'e Pierre et Marie Curie, Paris, France. \\
              \email{gregoire.nadin@upmc.fr} 
           \and 
            A.I. Toledo \at
              Laboratoire d'Analyse G\'eom\'etrie et Applications, Universit\'e Sorbonne Paris Nord, Villetaneuse, France.\\
              \email{anaisistoledo@gmail.com}           
           \and
           E. Ogier-Denis \at 
              Institut national de la santé et de la recherche médicale, Paris, France. \\
              \email{eric.ogier-denis@inserm.fr}
           \and
            H.Zaag \at
            Laboratoire d'Analyse G\'eom\'etrie et Applications, Universit\'e Sorbonne Paris Nord, Villetaneuse, France. \\
             \email{hatem.zaag@univ-paris13.fr}   
}

\date{Received: date / Accepted: date}

\maketitle

\begin{abstract}
    Crohn’s disease is an inflammatory bowel disease (IBD) that is not well understood. In particular, unlike other IBDs, the inflamed parts of the intestine compromise deep layers of the tissue  and are not continuous but separated and distributed through the whole gastrointestinal tract, displaying a patchy inflammatory pattern. In the present paper, we introduce a toy-model which might explain the appearance of such patterns. We consider a reaction-diffusion system involving bacteria and phagocyte and prove that, under certain conditions, this system might reproduce an activator-inhibitor dynamic leading to the occurrence of Turing-type instabilities. In other words, we prove the existence of stable stationary solutions that are spatially periodic and do not vanish in time. We also propose a set of parameters for which the system exhibits such phenomena and compare it with realistic parameters found in the literature. This is the first time, as far as we know, that a Turing pattern is investigated in inflammatory models.
\end{abstract}

\section{Introduction} \label{intro}

Ulcerative colitis and Crohn's disease represent the two main types of inflammatory bowel disease (IBD). Both are relapsing diseases and may present similar symptoms including long-term inflammation in the digestive system, however they are very different: Ulcerative colitis affects only the large intestine and the rectum whereas Crohn's disease can affect the entire gastrointestinal tract from the mouth to the anus. Typical presentations of Crohn's disease include the discontinuous involvement of various portions of the gastrointestinal tract and the  development  of  complications  including  strictures,  abscesses,  or  fistulas that compromise deep layers of the tissue while ulcerative colitis remains superficial but present no healthy areas between inflamed spots. 


There is consensus now that IBD result from an  unsuitable  response of a deficient mucosal immune system to the indigenous flora and other luminal antigens due to alterations of the epithelial barrier functions. We propose in this paper a simplified mathematical model aiming to recreate the immune response triggering inflammation. In the particular case of Crohn's disease, we seek to understand the patchy inflammatory patterns that differentiate patients suffering from this illness from those who has been diagnosed with ulcerative colitis.

IBD can be seen as an example of the acute inflammatory response of body tissues caused by harmful stimuli such as the presence of pathogenic germs or damaged cells. This protective response is also associated with the origin of other well-known diseases such as rheumatoid arthritis, the inflammatory phase in diabetic wounds or tissue inflammation, and has been extensively studied. Today it is still of central interest for researchers and, although several models have been proposed in order to understand the causes that lead to acute inflammation, the mathematical approach to this topic remains a recent field of research. A very complete review on the subject is provided in \cite{Vodovotz_2006,Vodovotz_2004}.

Among the mathematical works on inflammation we can refer to many models based on ordinary differential equation \cite{Day_2006,Dunster_2014,Herald_2010,Kumar_2004,Lauffenburger_1981,Mayer_1995,Reynolds_2006,Roy_2009,WENDELSDORF_2010}. Most of the authors take into account pro-inflammatory and anti-inflammatory mediators but also pathogens and other more or less realistic physiological variables. Depending on the parameters and the initial data these models manage to reproduce a variety of scenarios that can be observed experimentally and clinically; for example the case in which the host can eliminate the infection and also other situations in which the immune system cannot keep the disease under control or where the existence of oscillatory solutions determines a chronic cycle of inflammation. Most of the conclusions in the referenced papers are the result of stability study of the equilibrium states and numerical analysis of the simulations by phase portraits methods. In addition, in \cite{Day_2006,Kumar_2004,Roy_2009,WENDELSDORF_2010} a sensitivity analysis of the variables to the parameters of the models is performed in order to adjust the numerical results with experimental data and achieve greater biological fidelity of the model. 

Several authors had also considered spatial heterogeneity in order to model the inflammatory response, we can mention \cite{Khatib_2007,Khatib_2011,Ibragimov_2006} in the particular case of atherogenesis, \cite{Lauffenburger_1983,Penner_2012} in the tissue inflammation context and \cite{Chalmers_2015,Sullivan_2006} for the acute inflammatory response. The main variables of the models introduced in the mentioned works vary according to the dynamics that the authors wish to describe, the density of phagocytic cells, pro-inflammatory cytokines, anti-inflammatory mediators and bacteria are some standard quantities that are often taken into account. As in the ordinary differential equations approach the stability of the systems is systematically studied, in  \cite{Chalmers_2015,Consul_2014,Khatib_2007,Lauffenburger_1981} a vast analysis of all possible scenarios is performed depending on the values of the model parameters, the authors provide biological interpretation of such behavior as well as numerical simulations; furthermore, in \cite{Khatib_2011} the existence of travelling waves solutions is proved to be at the origin of a chronic inflammatory response. 

A different approach is presented in \cite{Penner_2012}, the model introduced in this paper aims to explain mathematically the patterns observed in the skin due to acute inflammation in the absence of specific pathogenic stimuli. By analyzing the stability of homogeneous and non-homogeneous states, sufficient conditions leading to the existence of such patterns solutions are obtained; several numerical examples are given as well. Similarly, in \cite{Lauffenburger_1983} authors claim that the instability of uniform steady distribution of phagocytic cells might trigger non-uniform cell density distributions which is potentially dangerous since tissue damages may occur in regions of high cell concentration. In this sense some sufficient conditions are given in order to prevent the existence of such kind of unstable states, these conditions primarily involve the phagocyte random motility coefficient and a chemotaxis coefficient included in the model.  

As suggested by \textit{in vitro} studies, phagocytic cells (big eaters) may move following a chemotactic impulse generated by the presence of pathogens germs, for this reason most of the authors cited above include the effect of chemotaxis by mean of the classical term first introduced by Patlak in 1953 and Keller and Segel in 1970 \cite{KellerSegel,Patlak_1953}. Nevertheless, there is no consensus on this assumption, as noted in \cite{Lauffenburger_1983}, \textit{in vivo}  observations more often show that the phagocytes seem to move within an infected lesion randomly, this is the case in the models introduced in \cite{Consul_2014,Khatib_2007,Khatib_2011}.

In the present paper, we propose a mechanism leading to patterns, which does not rely on chemotactism. We think the inflammatory response could be modeled by an activator-inhibitor system. Such systems are known to produce Turing mechanism, that is, periodic stationary solutions. This could possibly explain the patchy nature of Crohn's disease.

\section{The model}

We propose here a reaction-diffusion system modelling the dysfunctional immune response that triggers IBD. As mentioned in the introduction, this kind of systems have attracted much interest as a prototype model for pattern formation, in this case we refer in particular to inflammatory patterns. 

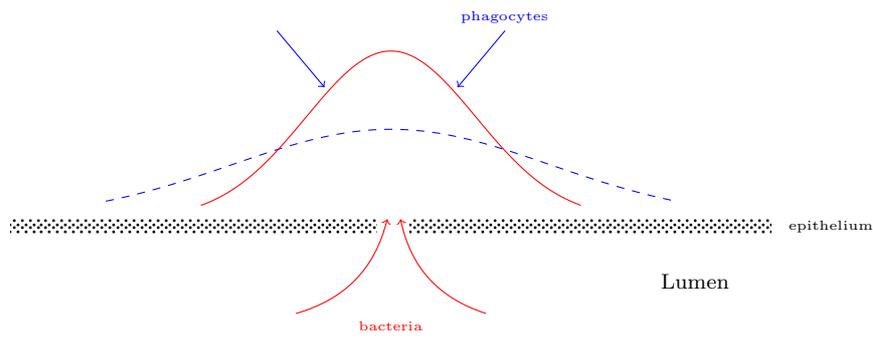
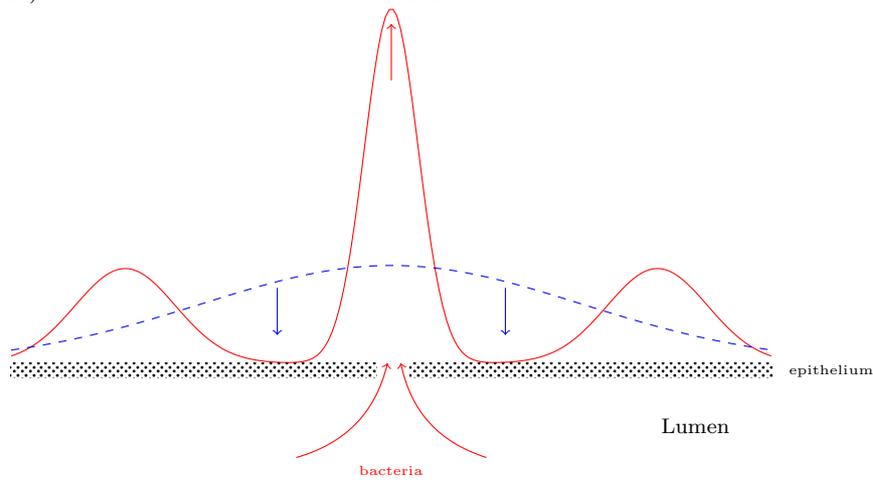
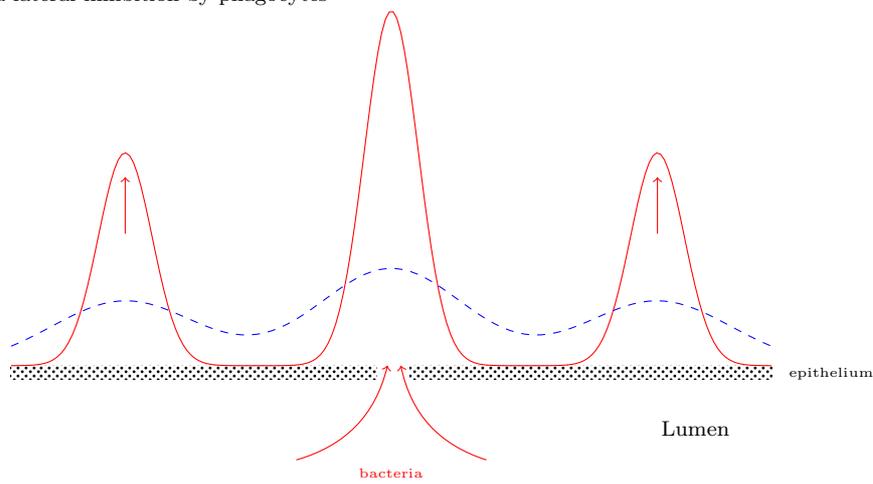
\begin{figure} [p]
	
	\begin{center}
		\begin{tikzpicture}[scale=2.5]
	
			\draw [domain=-1.5:1.5,blue,samples=100,dashed] plot(\x,{1/sqrt(0.7*2*3.14)*exp(-(\x)^2/(2*0.7))});
			\draw [domain=-1:1,red,samples=100] plot(\x,{1/sqrt(0.2*2*3.14)*exp(-(\x)^2/(2*0.2))});
			\fill [green, opacity=1,pattern=crosshatch dots] (-2,-0.075) rectangle (-0.075,0);
			\fill [green, opacity=1,pattern=crosshatch dots] (0.1,0) rectangle (2,-0.075);
			\draw [->,red] (0.5,-0.5) to[bend left,thick] (0.05,0);
			\draw [->,red] (-0.5,-0.5) to[bend right,thick] (-0.02,0);
			\draw (0,-0.5) node[below,fill=white] {\tiny \color{red} bacteria};
			\draw (0.6,1) node[above] {\tiny \color{blue} phagocytes};
			\draw (2.05,-0.04) node[right] {\tiny epithelium};
			\draw [->, blue] (0.6,1) to[thick] (0.35,0.7);
			\draw [->, blue] (-0.6,1) to[thick] (-0.35,0.7);
			\draw (1.6,-0.25) node[below] {Lumen};
		\end{tikzpicture} 
	\subcaption{Bacteria (red line) break through the epithelium (dotted zone); phagocytes (blue dashed line) are recruited in order to neutralize them}
    
	\begin{tikzpicture}[scale=2.5]
	
		\draw [domain=-2:2,red,samples=200] plot(\x,{1/(1.5*sqrt(0.02*2*3.14))*exp(-(\x)^2/(2*0.02))+1/(3*sqrt(0.07*2*3.14))*exp(-(\x-1.4)^2/(2*0.07))+1/(3*sqrt(0.07*2*3.14))*exp(-(\x+1.4)^2/(2*0.07))});
		\draw [domain=-2:2,blue,samples=100,dashed] plot(\x,{1.3/sqrt(2*3.14)*exp(-(\x)^2/2)});
		\fill [green, opacity=1,pattern=crosshatch dots] (-2,-0.075) rectangle (-0.075,0);
		\fill [green, opacity=1,pattern=crosshatch dots] (0.1,0) rectangle (2,-0.075);
		\draw [->,red] (0.5,-0.5) to[bend left,thick] (0.05,0);
		\draw [->,red] (-0.5,-0.5) to[bend right,thick] (-0.02,0);
		\draw (0,-0.5) node[below,fill=white] {\tiny \color{red} bacteria};
		\draw (2.05,-0.04) node[right] {\tiny epithelium};
		\draw [->, blue] (0.6,0.4) to[thick] (0.6,0.15);
		\draw [->, blue] (-0.6,0.4) to[thick] (-0.6,0.15);
		\draw [->, red] (0,1.5) to[thick] (0,1.8);
		\draw (1.6,-0.25) node[below] {Lumen};
	\end{tikzpicture} 
	\subcaption{Phagocytes spread rapidly through blood vessels. A high spot of bacteria remains with a lateral inhibition by phagocytes}

	\begin{tikzpicture}[scale=2.5]
	
		\draw [domain=-2:2,red,samples=200] plot(\x,{1/(1.5*sqrt(0.02*2*3.14))*exp(-(\x)^2/(2*0.02))+1/(2.5*sqrt(0.02*2*3.14))*exp(-(\x-1.4)^2/(2*0.02))+1/(2.5*sqrt(0.02*2*3.14))*exp(-(\x+1.4)^2/(2*0.02))});
		\draw [domain=-2:2,blue,samples=100,dashed] plot(\x,{1/(2*sqrt(0.15*2*3.14))*exp(-(\x)^2/(2*0.15))+1/(3*sqrt(0.15*2*3.14))*exp(-(\x-1.4)^2/(2*0.15))+1/(3*sqrt(0.15*2*3.14))*exp(-(\x+1.4)^2/(2*0.15))});
		\fill [green, opacity=1,pattern=crosshatch dots] (-2,-0.075) rectangle (-0.075,0);
		\fill [green, opacity=1,pattern=crosshatch dots] (0.1,0) rectangle (2,-0.075);
		\draw [->,red] (0.5,-0.5) to[bend left,thick] (0.05,0);
		\draw [->,red] (-0.5,-0.5) to[bend right,thick] (-0.02,0);
		\draw (0,-0.5) node[below,fill=white] {\tiny \color{red} bacteria};
		\draw (2.05,-0.04) node[right] {\tiny epithelium};
		\draw [->, red] (1.4,0.7) to[thick] (1.4,1);
		\draw [->, red] (-1.4,0.7) to[thick] (-1.4,1);		
		\draw (1.6,-0.25) node[below] {Lumen};
	\end{tikzpicture} 
	\subcaption{Other spots appear}
	
\end{center}
\caption {Initiation of the inflammatory process} \label{scheme}
\end{figure}

Roughly speaking, the first line of defense of the mucosal immune system is the epithelial barrier which is a polarized single layer covered by mucus in which commensal microbes are embedded. Lowered epithelial resistance and increased permeability of the inflamed and non-inflamed mucosa  is systematically observed  in  patients  with  Crohn's  disease and ulcerative colitis, hence the epithelial   barrier gets leaky and luminal antigens gain access to the underlying mucosal tissue. In a healthy gut, the immune response by mean of intestinal phagocytes  eliminates the external agents limiting the inflammatory response in the gut. Unfortunately in a disease-state the well controlled balance of the intestinal immune system is disturbed at all levels, this dysfunctional mechanism contributes to acute and chronic inflammatory processes. Indeed, an excessive amount of immune cells migrating to the damaged zone can engage the permeability of the epithelial barrier and thus might allows further infiltration of microbiota which aggravate inflammation. This complex network triggers the initiation of an inflammatory cascade that causes ulcerative colitis and Crohn's diseases, see Fig.\ref{scheme}.

For the sake of simplicity in this model we will consider just two components varying in time and space: {\footnotesize 1.} The number of non-resident bacteria leaking into the intestinal tissue through the epithelial barrier noted as $\beta$, also refereed as microbiota, pathogens or antigens and  {\footnotesize 2.} The immune cells $\gamma$ which we will often refer as phagocytic cells. Also, by simplicity we model a portion of the digestive tube as an interval $\Omega\subset \R$ of the real axis, which will be very large. The model reads:
\begin{equation}\label{Modelo}\left\{
	\begin{array}{lcl}
		\partial_{t} \beta - d_{b}\Delta \beta &=& r_{b}\left(1-\frac{\beta}{b_i} \right)\beta -\frac{a \beta \gamma}{s_b+\beta} +f_e \left(1-\frac{\beta}{b_i}\right)\gamma,\\
		\partial_{t} \gamma - d_{c}\Delta \gamma &=& f_{b} \beta -r_c \gamma.\\
	\end{array}\right.
\end{equation}
We complete by considering  Neumann boundary conditions and initial data $\beta(0,x)=\beta_0(x)$ and $\gamma(0,x)=\gamma_0(x)$ for all $x \in \Omega$.

During the immune response there is a first stage where the non-resident phagocytes migrate from the vasculature into  the  intestinal  mucosa and a second stage where they move to the damaged zone and fight the bacteria. This first stage results from a transport movement through the blood vessels and it is almost instantaneous compared to the second one, so we omit it in this simplified model. 

Another main assumption is to consider that immune cells and bacteria move randomly through the damaged tissue and the epithelial barrier. As mentioned in the introduction, it is generally accepted that diffusion provides an adequate description of molecular spreading but, in the case of phagocytic cells, chemotaxis is claimed to be crucial establishing the direction of  movement in the sense of the pathogen gradient. However, there are \textit{in vivo} experiments that corroborate our hypothesis \cite{Lauffenburger_1983} and several authors have made similar assumptions \cite{Consul_2014,Khatib_2007,Khatib_2011}. Nevertheless, by neglecting chemotaxis in our model we do not claim that it is an unimportant phenomenon, instead, this assumption must be seen as a simplification and an idealization of the physiological mechanism we seek to describe.  

The coefficients $d_b>0$ and $d_c>0$ are the diffusion rates of bacteria and phagocytes, respectively. The parameter $r_b>0$ is associated with the reproduction rate of bacteria. 

In healthy conditions the number of bacteria within the lumen remains almost constant and they are not able to penetrate the epithelial barrier, we associate this quantity to the parameter $b_i>0$.  We remark that this parameter $b_i$ is in some sense a carrying capacity; in fact, in the total absence of the epithelial barrier, the maximum amount of bacteria in the colon would not be greater than $\beta=b_i$, that is the reason why we add the logistic term $1-\frac{\beta}{b_i}$ in the first equation, \cite{Verhulst_1845,Perthame_2015}. 

The parameter $f_b>0$ is associated with the immune response rate of the organism sending cells to fight bacteria in the damaged zones. In others words, as soon as the presence of pathogens is detected, phagocytes are coming up.  

The term $-\frac{a \beta \gamma}{s_b+\beta}$ with $a>0$ and $s_b>0$ corresponds to the effect of the immune system on the pathogen agents. In particular $\frac{a \beta }{s_b+\beta}$ is the phagocytosis rate or intake rate, it suggests that the attack rate of immune cells on bacteria varies with the density of pathogen. This functional response term takes into account the rate $p_{c}$ at which phagocytes encounter a bacterium per unit of bacteria density, which is $p_{c}:=\frac{a}{s_b}$ and  the average time $\tau$ that takes a phagocyte to neutralize a bacterium (or handling time) which can be computed as $\tau:=\frac 1{a}$. Experiments presented in \cite{Leijh_1980,Stossel_1973} reflect this dynamic. In the mathematical literature such kind of term is often referred as a Holling Type II functional response, see  \cite{holling_1965,Perthame_2015}. 

We consider $f_e>0$ as a measure of the negative effect of the phagocyte's concentration for the epithelial resistance, and therefore it has a positive impact on the bacteria density i.e. the larger the epithelial gap, the more there are bacteria, the more there are immune cells drifting to the damaged zone and the more porous is the epithelium and so on. 

Finally, a self-regulation function of anti-inflammatory cells limits their life-time, so immune cells have an intrinsic death rate which is noted in the model as $r_c>0$.

\section{On Turing Patterns}

Since one of our main interest with this paper is to explain patchy inflammatory bowel patterns often observed in patients suffering from Crohn's disease, we seek to demonstrate that the model we propose may present Turing-type instabilities under certain conditions. This denomination is due to Alan Turing who was the first to describe spatial patterns caused by the effects of diffusion in his article on morphogenesis theory published in 1952, \cite{Turing_1952}.

Roughly speaking, a Turing system consist of an activator that must diffuse at a much slower rate than an inhibitor to produce a pattern. We remind to the reader that diffusion causes areas of high concentration to spread out to areas of low concentration. In such kinds of systems the activator component must increase the production of itself while the inhibitor restrains the production of both. Turing's analysis shows that in certain regimes those systems are unstable to small perturbations, leading to the growth of large scale patterns.

In the model we previously introduce bacteria are the activator and the immune cells the inhibitor, indeed bacteria reproduce at a certain rate $r_b$ and immune cells neutralize bacteria by phagocytosis (Holling-type term) and self-regulate their own life-time $r_c$. In practice, we should look for steady state solutions of the equation (\ref{Modelo}) which are linearly unstable, i.e. such that there are perturbations for which the linearized system has exponentially growing solutions in time. To be sure that a Turing-type phenomena is occurring it is important to exclude the cases where the corresponding growth modes are unbounded, that is solutions with infinitely high frequencies and also the cases in which solutions blow up or go to extinction \cite{Perthame_2015}. 

In section \ref{Stability and Turing patterns} we study the conditions leading to the observation of Turing phenomena in our model.

\section{Results}

\subsection{Non-negativity property and boundedness}

We begin by establishing some elementary properties in the model to guarantee system (\ref{Modelo}) accuracy as a population dynamics model. In other words it is important that whenever the initial data have a reasonable biological meaning, the solution of the differential equation inherits that property. We start by a non-negativity property:

\begin{prop}\label{positivite}
	Provided that the initial condition $(\beta_0(x),\gamma_0(x))$ is non-negative the solutions of the system (\ref{Modelo}) remain non-negative for every $t>0$.
\end{prop}

Similarly, we establish a boundedness property associated with the carrying capacity of the population environment:

\begin{prop}\label{bornes}
	If $\beta_0(x)<b_i$ then for all $t>0$ one has $\beta(t,x)<b_i$. Moreover, if $\gamma_0(x)$ is bounded, then $\gamma(t,\cdot)$ remains bounded in the $L^2$-norm in $\Omega$ for every $t>0$.
\end{prop}

\subsection{Stability analysis} \label{Stability and Turing patterns}
Let us study now the steady states of the model and their stability properties. The equation (\ref{Modelo}) have two non-negative homogeneous steady states. One of them is the trivial solution $(\beta,\gamma)=(0,0)$ associated with the absence of bacteria and immune cells. The other one, that we denote  $(\beta,\gamma)=(\overline{\beta},\overline{\gamma})$, satisfies:

\begin{equation}\label{caracterizacion}
0=\left(r_{b}+f_e \kappa \right)\left(1-\frac{\overline{\beta}}{b_i}\right) -\frac{a \kappa \overline{\beta}}{s_b+\overline{\beta}},
\end{equation}
where $\kappa:=\frac{f_b}{r_c}$ and $\overline{\gamma}=\kappa\overline{\beta}$.

We remark that $(0,0)$ is unstable. Indeed, the linearized matrix around this steady state has negative determinant and thus an eigenvalue with positive real part. For the non-trivial equilibrium point $(\overline{\beta},\overline{\gamma})$ the stability analysis is less straightforward. The following proposition establishes the conditions leading to the stability of this steady state.

\begin{prop}\label{ODE_linear_stability}
	Consider the O.D.E system associated with (\ref{Modelo}) with non-negative real parameters $a, r_b, r_c, f_b, f_e, b_i$ and $s_b$,
	\begin{equation}\label{Modelo_ODE}\left\{
	\begin{array}{lcl}
	\displaystyle{\partial_{t} \beta} &=&  \displaystyle{r_{b}\left(1-\frac{\beta}{b_i}\right)\beta -\frac{a \beta \gamma}{s_b+\beta} +f_e \left(1-\frac{\beta}{b_i}\right)\gamma}\\
	\displaystyle{\partial_{t} \gamma } &=&\displaystyle{ f_{b} \beta -r_c \gamma}\\
	\end{array}\right..
	\end{equation}
	This system has a unique positive steady state solution  $\big( \beta(t),\gamma(t)\big)=(\overline{\beta},\overline{\gamma})$ which is stable if and only if 
	\begin{equation} \label{C2}
	\frac{a\kappa \overline{\beta}^2}{(s_b+\overline{\beta})^2} - r_b \frac{\overline{\beta}}{b_i} - f_e\kappa < r_c.
	\end{equation}
\end{prop}

We conjecture that the model might show some unexpected behavior around this steady state which could be at the origin of patchy inflammatory patterns. Hence, let us focus on conditions leading the formation of Turing patterns for the reaction diffusion system (\ref{Modelo}), that is perturbations around the steady state $(\overline{\beta},\overline{\gamma})$ such that the linearized system has exponential growth in time and for which the corresponding growth modes are bounded. The following proposition establishes the necessary conditions for the occurrence of such phenomenon. 

\begin{prop}\label{Turing_patterns}
	Consider the system (\ref{Modelo}) and its unique positive homogeneous steady state solution $(\overline{\beta},\overline{\gamma})$; assume that there exist real non-negative values of the parameters $a, r_b, r_c, s_b, f_e, f_b, b_i$ such that the following condition holds:
	\begin{equation}\label{Conditions}
	0 < \frac{a\kappa \overline{\beta}^2}{(s_b+\overline{\beta})^2} - r_b \theta - f_e\kappa  < r_c
	\end{equation}
	Then for $\frac{d_b}{d_c}$ small enough the reaction diffusion system (\ref{Modelo}) shows Turing instabilities around this steady state.
	
\end{prop}

\section{Parameters of the model}

In this section we want to estimate the values of the parameters of the model and to prove the non emptiness of the parameters set defined by (\ref{Conditions}). As long as it is possible we will rely on values obtained from real observations or in vitro experiments. However, in some cases the exact values are unknown due to the difficulty of measuring them in vivo or even in vitro.

Let us start with an estimation of the reproduction rate of the bacteria, represented in our model as $r_b$. Bacterium's generation time, which is the time it gets to the population to double the number of individuals, might vary from 12 minutes to several hours depending on temperature, nutrients, culture medium, among others factors. For E. Coli, for instance,  it is around 20 minutes in standard conditions, \cite{Korem_2015}. We can then consider that  the evolution of bacteria population is given by $\partial_t b = r_b b$ and so $r_b = \frac{ln(2)}{20}$ measured in bacteria per minute. That gives us an approximate value $r_b = 3.47\ast 10^{-2}$ $u/min$ which is in the estimated range of values given in \cite{Lauffenburger_1983} for this parameter.

Similarly, it is known that in healthy conditions phagocytes have, in average, a half-life of two days \cite{Labro_2000}, and so from $\partial_t c = -r_c c$ we get $r_c = \frac{ln(2)}{2880}$ cells per minute which means that the death rate of phagocytes is ideally of the order of $10^{-4}$ $u/min$, which coincides with that considered in \cite{Waugh2007} for immune cells in diabetic wounds or in \cite{Lauffenburger_1983} for bacterial infection causing tissue inflammation. However, there is no consensus, some authors assume this parameter to be of the order of $10^{-3}$$u/min$ in the inflammatory response framework \cite{Chow_2005} or even of the order of $10^{-6}$$u/min$ in the case of early atherosclerosis \cite{Chalmers_2015}. For such parameters, corresponding to a healthy organism, we do not expect to observe a Crohn's disease. Indeed, the mechanism we describe below occurs with $r_c=2\times 10^{-2}$ $u/min$ (see Table \ref{parameters}). For $r_c=10^{-3}$ $u/min$, the range of parameters for which a Turing pattern occurs is quite narrow, Fig.\ref{figura3}.   

The diffusion coefficient of immune cells might also vary according to the type of cell and the part of the body where they act. In the consulted literature the value of this parameter varies from $10^{-12}$ $m^2/min$ to $10^{-10}$ $m^2/min$ depending on the context \cite{Consul_2014,Khatib_2011,Lauffenburger_1983,Stickle1985}. In the absence of experimental data providing more precise information about the order of this parameter in the particular case of bacterial infection in the intestinal track, we consider this coefficient to remain within this range in damaged areas of the intestine. 

Although there are not precise information concerning the diffusion rate of bacteria through the epithelial barrier, it is known that in aqueous solutions like the lumen, the diffusion rate might vary from $10^{-11} $ $m^2/min$ to $10^{-8} $ $m^2/min$ depending on the type of bacteria. However, in a non-liquid  framework, which is the case of bacteria penetrating through the epithelial barrier, motility should be reduced.  

We will now roughly compute a value for the parameter $a$, we suppose that there is a signficant density of bacteria in a certain position $x=x_0$, and we study the time evolution of the population within this point. If $\beta$ is large enough, the term $1-\frac{\beta}{b_1}$ is negligible, moreover the term $-\frac{a\beta\gamma}{s_b+\beta}$ tends to approach $-a \gamma$, so we can approximately write 
\begin{equation} \label{approximacion_de_alpha}
   \partial_{t}\beta(t,x_0)=-a \gamma(t,x_0).
\end{equation}
Let us now define $\tau$ as the average time it takes a phagocyte to neutralize a bacterium, which is around 3 minutes in the in vitro observations, it implies that 
\begin{equation}
	\beta(t+\tau,x_0)=\beta(t,x_0)-\gamma(t,x_0)
\end{equation}
and consequently $\partial_{t}\beta(t,x_0)\approx -\frac{\gamma(t,x_0)}{\tau}$. Replacing this into (\ref{approximacion_de_alpha}) we arrive at the conclusion that $a$ is of the order of $\frac{1}{\tau}$ units per minute.

The density of bacteria in the lumen is approximately $b_i = 10^{17}$ $u/{m^3}$. At the positive equilibrium stage $(\overline{\beta},\overline{\gamma})$,which is associated to an inflammatory phase, we suppose that around 30$\%$ of the total density of bacteria within the lumen might penetrate the epithelial barrier without going out of control. Therefore, we set $\overline{\beta}=0.3\times b_i $ units of bacteria. Even though we have no exact data concerning the density of immune cells in the damaged zone, the in vitro experiments suggest that during the inflammation stage it is around ten times less than the bacteria density, this is quite natural considering that the size of a phagocyte is much larger than the size of a bacterium. Hence, we set the hypothesis that $\kappa=\frac{1}{10}$ which means that at the equilibrium point it holds $\overline{\beta} = 10\overline{\gamma}$ and consequently $\overline{\gamma}= 3\cdot 10^{-2} \times b_i$. Taken this into account from the equilibrium condition we have that $f_b = 10^{-1} r_c$ measured in units per minute. 

The parameter $f_e$ is finally computed so that (\ref{caracterizacion}) holds at the equilibrium state. 

\section{Numerical simulations}

\begin{figure}[t]
	\includegraphics[width=0.95\textwidth]{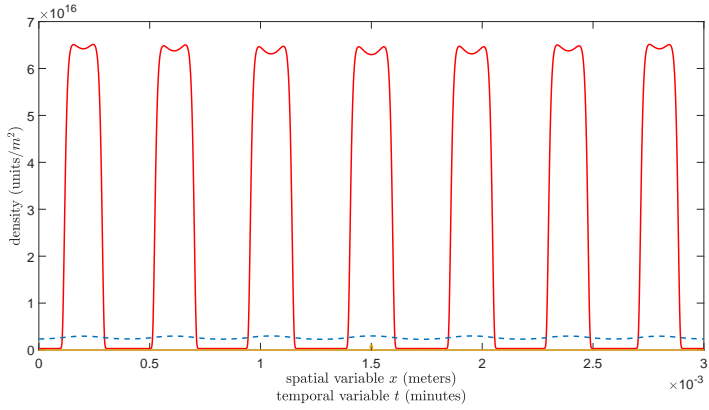}
	\caption{Bacteria (red line) and phagocytes (blue dashed line) after a time-lapse of 2 weeks with an initial bacterial infection  $\beta_0(x)=10^{15}\times\mathbb{I}_{\{1495 \leq x \leq 1505\}}$ (yellow line) and $\gamma_0(x)=0$}\label{figura2}
\end{figure}

\begin{figure}
	\includegraphics[width=0.95\textwidth]{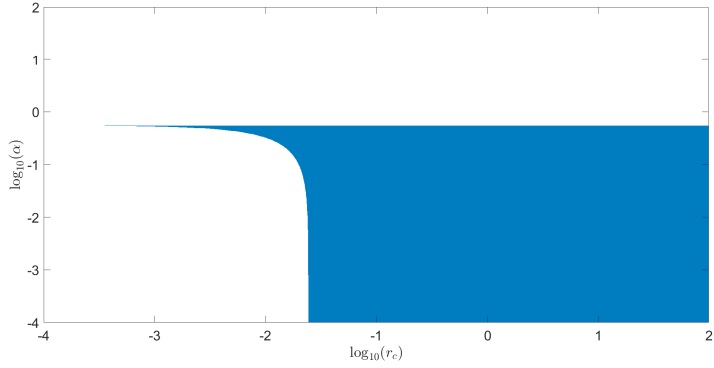}
	\caption{Region (blue) defined by the parameters $r_c$ and $a$ that verify condition (\ref{Conditions}) leading to Turing patterns observation} \label{figura3}
\end{figure}

We perform some numerical simulations in MATLAB by mean of a semi-implicit scheme to solve the system of equations (\ref{Modelo}), the results are shown in Fig. \ref{figura2}. We have considered the parameters values presented in the table \ref{parameters} which were estimated in the previous section. For these values, the condition (\ref{Conditions}) associated to a Turing phenomenon occurrence established in the Proposition \ref{Turing_patterns} is verified. However, there is a whole family of parameters verifying (\ref{Conditions}), as shown in Fig. \ref{figura3}.

\begin{table}[b]
    \caption{Assigned values for the parameters of the model (\ref{Modelo})}
    \label{parameters}  
    \centering    
    \begin{tabular}{llll}
    \hline\noalign{\smallskip}
    parameter & interpretation & value & units \\
    \noalign{\smallskip}\hline\noalign{\smallskip}
    $r_b$  & Reproduction rate of bacteria & 0.0347 & \tiny{(u/min)}\\
    $r_c$ &  Intrinsic death rate of phagocytes & 0.02 & \tiny{(u/min)} \\
    $d_b$ & Diffusion rate of bacteria & $10^{-13}$ & \tiny{($m^2$/min)}\\
    $d_c$  & Diffusion rate of phagocytes & $10^{-10}$ & \tiny{($m^2$/min)}\\
    $b_i$  & Density of bacteria in the lumen & $10^{17}$ & \tiny{(u/$m^3$)} \\
    $f_b$  & Immune response rate & 0.002 & \tiny{(u/min)}\\
	$a$  & Coefficient proportional to the rate of phagocytosis ($a=s_b p_{c}$) \ \  & 0.3129 & \tiny{(u/min)} \\
	     &  it is also inversely proportional to the handling time ($a=\frac 1{\tau}$) & &  \\
    $s_b$  & Proportionality coefficient between $p_{c}$ and $a$  & $10^{15}$ & \tiny{(u/$m^3$)}\\
    $f_e$  & Related to the porosity of the epithelium & 0.0856 & \tiny{(u/min)} \\
    \noalign{\smallskip}\hline
    \end{tabular}
\end{table}

For the simulations we have considered an initial datum with no phagocytes presence and a tiny spot of bacteria concentrated in the middle of the domain $\Omega$. This might be understood as a slight leak of bacteria from the lumen through the epithelium. The activator-inhibitor dynamics generated by the body's immune response to the presence of bacteria and the contrast in the propagation rates of the two actors of the system is the reason why the patterns emerge in Fig. \ref{figura2} after a certain time. This behavior is definitively associated with a Turing phenomenon.

We remark that the values we assign to the diffusion coefficients remains within the range estimated in the previous section. However, from the mathematical point of view what is really important in order to ensure verifying the conditions leading to the observation of Turing patterns is the smallness of the ratio $\delta = \frac{d_b}{d_c}$. To change those values by preserving $\delta$ only represents a spatial rescaling that does not affect the pattern formation.

\section{Proof of the results}

\subsection*{Proof of the Proposition \ref{positivite}}
\begin{proof}
	Consider $\overline{t}>0$ the first instant when either $\beta(t,x)$ or $\gamma(t,x)$ became non-positive, then for some $x^*\in \Omega$ one has $\beta(\overline{t},x^*) \gamma(\overline{t},x^*)=0$. 
	
	If $\beta(\overline{t},x^*)=0$ and $\gamma(\overline{t},x^*)\geq 0$ since $\beta_0(x)>0$ then there exist $\delta_t>0$ such that 
	\begin{equation}
	\forall t\in [\overline{t}-\delta_t,\overline{t}] \ \hbox{one has} \ \partial_{t}\beta(t,x^*)<0.
	\end{equation}
	Nevertheless, from the first equation in (\ref{Modelo}) one has $\partial_{t}\beta(\overline{t},x^*)=f_e \gamma(\overline{t},x^*)\geq 0$ which contradicts the previous conclusion.
	
	Similarly, if $\beta(\overline{t},x^*)\geq 0$ and $\gamma(\overline{t},x^*)= 0$ from the positivity assumption of the initial data we can conclude the existence of $\delta_t>0$ such that 
	\begin{equation}
	\forall t\in [\overline{t}-\delta_t,\overline{t}] \ \hbox{one has} \ \partial_{t}\gamma(t,x^*)<0. 
	\end{equation}
	Again from the second equation in (\ref{Modelo}) one has $\partial_{t}\gamma(\overline{t},x^*)=f_b \beta(\overline{t},x^*)\geq 0$ which is a contradiction.
\end{proof}

\subsection*{Proof of the Proposition \ref{bornes}}
\begin{proof}
	The argument of this proof is similar to the one used to prove the non-negativity property. Indeed, consider $\overline{t}$ the first instant when $\beta$ rises the value $b_i$, then there exist $x^*\in \Omega$ such that $\partial_t \beta(\overline{t},x^*)>0$, nevertheless from the equation associated with $\beta$ one conclude that $\partial_t \beta(\overline{t},x^*)=-\frac{a \beta(\overline{t},x^*)\gamma(\overline{t},x^*)}{s_b+\beta(\overline{t},x^*)}<0$ from the positivity property. So we get a contradiction which implies that for all $t>0$ one has necessarily $\beta(t,x)<b_i$.
	
	The boundedness of $\gamma$ follows directly from the boundedness of $\beta$ and $\gamma_0$. In fact multiplying by $\gamma$ in the second equation of (\ref{Modelo}), integrating by parts and applying Holder inequality one gets that
	\begin{equation}
	\frac{1}{2}\partial_{t} \|\gamma\|^2_{L^2(\Omega)}+\|\nabla\gamma\|^2_{L^2(\Omega)}\leq \left( f_b\|\beta\|^2_{L^2(\Omega)}-r_c\right)\|\gamma\|^2_{L^2(\Omega)}
	\end{equation}
	from where after applying the Gronwall inequality one concludes that there exist a positive constant $C=C(t)$ such that  $\|\gamma(t)\|^2_{L^2(\Omega)}<C\|\gamma_0\|^2_{L^2(\Omega)}$.
\end{proof}

\subsection*{Proof of the Proposition  \ref{ODE_linear_stability}}
\begin{proof} 
	The existence of such a positive steady state follows from the analysis of (\ref{caracterizacion}). Let us define $F(\beta)=\left(r_{b}+f_e \kappa \right)\left(1-\frac{\beta}{b_i}\right) -\frac{a \kappa \beta}{s_b+\beta}$. From the positivity of the parameters of the model we have that $F(0)>0$ and $F(b_i)<0$, this means that there are at least one positive value $\overline{\beta}\in(0,b_i)$ that satisfies $F(\overline{\beta})=0$ or equivalently (\ref{caracterizacion}). Moreover, since the derivative of $F$ is strictly negative we deduce that it has at most one root which leads to the uniqueness of $\overline{\beta}$.
	
	Let us now study the conditions leading to the stability of this steady state. In order to simplify the notations we will define $\theta:=\frac{\overline{\beta}}{b_i}$. We also define $\mathbf{M}$ as the matrix of the linearized system around this positive steady state  $(\overline{\beta},\overline{\gamma})$ 
	
	\[
        \mathbf{M}:=
	\begin{pmatrix}
	\displaystyle{r_b(1-2\theta)-\frac{a s_b \kappa \overline{\beta}}{(s_b+ \overline{\beta})^2}-f_e \kappa \theta} \  & \  \displaystyle{-\frac{a \overline{\beta}}{s_b + \overline{\beta}} + f_e(1-\theta)} \\
	\displaystyle{f_b} & -r_c
	\end{pmatrix}.
	\]
	We compute the determinant and the trace of this matrix
	\begin{eqnarray*}
		tr(\mathbf{M}) &=& \frac{a\kappa \overline{\beta}^2}{(s_b+\overline{\beta})^2} - r_b \theta - f_e\kappa  - r_c,\\
		det(\mathbf{M})&=& r_cr_b \theta + f_bf_e\theta + \frac{af_bs_b \overline{\beta}}{(s_b + \overline{\beta})^2}.\\
	\end{eqnarray*}
	From the positivity of the parameters of the model it is clear that the determinant of $\mathbf{M}$ is positive, therefore in order to have linear stability around $(\overline{\beta},\overline{\gamma})$ it is necessary and sufficient to impose the negativity of the trace of $\mathbf{M}$ which is equivalent to (\ref{C2}).
	
\end{proof}

\subsection*{Proof of the Proposition  \ref{Turing_patterns}}

\begin{proof}
	We linearize the system around $(\overline{\beta},\overline{\gamma})$. For the sake of simplicity we keep the notation $\beta(t,x),\gamma(t,x)$ for the linearized variables
	\begin{equation}\left\{
	\begin{array}{lcl}
	\partial_t \beta - d_b \Delta \beta &=& \bigg(r_b(1-2\theta)-\frac{a\kappa s_b \overline{\beta}}{(s_b + \overline{\beta})^2}-f_e \kappa \theta \bigg) \beta + \bigg(-\frac{a \overline{\beta}}{s_b + \overline{\beta}} + f_e(1-\theta)\bigg)\gamma \\
	\partial_t \gamma - d_c \Delta \gamma &=& f_b \beta -r_c \gamma
	\end{array}\right.
	\end{equation}
	
	We are seeking in particular for solutions with exponential growth in time, so we consider that 
	\begin{equation}\label{descomposicion}
	\beta(t,x)=e^{\lambda t}B(x) \ ; \  \gamma(t,x)=e^{\lambda t}C(x)
	\end{equation}
	with $\lambda > 0$. This means that $B(x)$ and $C(x)$ should satisfy the fallowing problem  
	\begin{equation}\label{linear_system}\left\{
	\begin{array}{lcl}
	- d_b \Delta B(x) &=& \bigg(r_b(1-2\theta)-\frac{a\kappa s_b \overline{\beta}}{(s_b + \overline{\beta})^2}-f_e \kappa \theta - \lambda \bigg) B(x) + \bigg(-\frac{a \overline{\beta}}{s_b + \overline{\beta}} + f_e(1-\theta)\bigg) C(x) \\
	- d_c \Delta C(x) &=& f_b B(x) +(-r_c -\lambda) C(x)
	\end{array}\right.
	\end{equation}
	or equivalently that they are eigenfunctions associated with the positive eigenvalue $\lambda$. 
	We consider in particular Fourier modes of the form 
	\[B(x)=Be^{i\xi x} \quad ; \ C(x)=Ce^{i \xi x}, \]
	and we replace it in (\ref{linear_system}) to obtain the fallowing homogeneous linear system of equations  
	\[
	\begin{pmatrix}
	0 \\
	0 \end{pmatrix}= 
	\begin{pmatrix}
	r_b(1-2\theta)-\frac{a\kappa s_b \overline{\beta}}{(s_b + \overline{\beta})^2}-f_e \kappa \theta - \lambda -d_b \xi^2 \  & \ -\frac{a \overline{\beta}}{s_b + \overline{\beta}} + f_e(1-\theta)  \\
	f_b &  -r_c - \lambda - d_c \xi^2
	\end{pmatrix}\begin{pmatrix}
	B \\
	C \end{pmatrix}
	\]

	Let us call $\mathbf{M}_{\lambda,\xi}$ the matrix associated  to the previous linear system. It can be written in terms of $\xi$, $\lambda$ and the matrix $\mathbf{M}$ introduced before in the proof of the Proposition \ref{ODE_linear_stability}
	\[
        \mathbf{M}_{\lambda,\xi}= 
	\begin{pmatrix}
	\displaystyle{\mathbf{M}_{(1,1)}  - \lambda - d_b\xi^2} \  & \  \mathbf{M}_{(1,2)} \\
	\mathbf{M}_{(2,1)} &  \mathbf{M}_{(2,2)} - \lambda -d_c \xi^2
	\end{pmatrix}.
	\] 
	In other words we look for a certain $\lambda$ with positive real part and $\xi^2$ for which $\det(\mathbf{M}_{\lambda,\xi})=0$. The determinant of $\mathbf{M}_{\lambda,\xi}$ is a quadratic polynomial function in $\lambda$
	\begin{equation}\label{determinante}
	det(\mathbf{M}_{\lambda,\xi}) = \lambda^2 + a_1\lambda + a_2
	\end{equation}
	with coefficients 
	\begin{eqnarray*}
		a_1 &=& -tr(\mathbf{M}) +(d_b+d_c)\xi^2,\\
		a_2 &=& \det(\mathbf{M}) - (\mathbf{M}_{(1,1)}d_c+\mathbf{M}_{(2,2)}d_b)\xi^2+d_bd_c\xi^4.
	\end{eqnarray*}
	Since the right-hand side inequality in (\ref{Conditions}) ensures that $tr(\mathbf{M})<0$, we conclude that $a_1>0$. Hence, the polynomial associated to $\det(\mathbf{M}_{\lambda,\xi})$ can have a positive root $\lambda$ if and only if $a_2<0$. The term $a_2$ is itself a quadratic polynomial in $\xi^2$ with positive second order coefficient. For the sake of simplicity we will define $\delta:=\frac{d_b}{d_c}$, and we will study the sign of $\frac{a_2}{d_bd_c}$ which roots are explicitly given by 
	\begin{equation}
	\Lambda_{\pm}=\frac{\mathbf{M}_{(1,1)}+\delta \mathbf{M}_{(2,2)}}{2*d_b}\bigg[1 \pm \sqrt{1-\frac{4\det(\mathbf{M}) \delta}{(\mathbf{M}_{(1,1)}+\delta \mathbf{M}_{(2,2)})^2}} \bigg].
	\end{equation}
	In the regime $\delta$ small enough the Taylor expansion gives us the following approximate values 
	\begin{eqnarray}
	\Lambda_-&=&-\frac{\det(\mathbf{M})}{d_c \mathbf{M}_{(1,1)}}\\
	\Lambda_+&=&\frac{\mathbf{M}_{(1,1)}}{d_c \delta}
	\end{eqnarray} 
    The left-hand side inequality in (\ref{Conditions}) guarantees that $\Lambda_+$ is positive and since $\delta$ can be as small as desired, then $\Lambda_+>>1$ and the interval $(\Lambda_-,\Lambda_+)$ where $a_2$ is negative is large enough.
    
	In other words,  there exist a positive real $\lambda$ and Fourier modes for which $\det(\mathbf{M}_{\lambda, \xi^2})=0$ and consequently we can find exponential growth in time solutions to the linearized system around the steady state $(\overline{\beta},\overline{\gamma})$. However, the Fourier modes for which this condition holds are bounded. 
	
	We have showed the existence of perturbations such that the linearized system has exponential growth in time. The frequency of the perturbations can not be infinity and from Proposition \ref{positivite} and \ref{bornes} neither extinction nor blows-up are possible. Hence, we have finally proved the formation of Turing Patterns.
	
\end{proof}

\section{Conclusions}
 This work remains a simplified approach to the question of modelling inflammatory response in Crohn's disease. We have made several hypotheses with the aim of globally understanding the biological mechanism behind the abnormal body reaction leading to the disease but staying relatively simple in terms of number of variables and equations. 

 Though we have tried to consider parameters values true to medical and biological observations, we highlight the qualitative results over quantitative ones. In this sense, obtaining a Turing mechanism through our model, might explain the patchy inflammatory patterns observed in patients  suffering from Crohn's disease and must be interpreted as another step in the aiming to fully understand this illness and its causes.

 It remains a question concerning the Ulcerative Colitis (RCH) since it has several common factors that relate it to Crohn's disease but also others that set them apart. It might be interesting to study the possibility of modelling RCH by mean of the same system of equations in a different parameter regime and eventually find responses helping doctors with early diagnosis or treatments.


\begin{acknowledgements}
{This work would not have been possible without the support of the Inflamex Laboratory of Excellence and the Galilee PhD College, whom we sincerely thank. Also, a special thanks to Dr. Xavier Treton for helping us understand inflammatory bowel diseases.}
\end{acknowledgements}

%

\paragraph{Declarations}

\section*{Funding}

{This work was supported by the Inflamex Laboratory of Excellence.}

\section*{Conflict of interest}

The authors declare that they have no conflict of interest.

\bibliographystyle{apalike}
\bibliography{paper_Crohn_Arxiv}

\end{document}